# A CONFORMAL CO-SYMPLECTIC STRUCTURE ON THE SPACE OF PSEUDO-RIEMANNIAN GEODESICS

PATRICK IGLESIAS-ZEMMOUR

ABSTRACT. The classical construction of the symplectic structure on the space of geodesic trajectories via Hamiltonian reduction fails in the pseudo-Riemannian setting due to a dimensional mismatch created by the null geodesics. This paper proposes a new, unified approach. We first construct the space of all geodesic trajectories $\mathscr{G}_{\text{traj}}$ directly as the quotient of the space of geodesics curves $\mathscr{G}_{\text{curv}}$ by the affine reparametrization group. We then perform a symplectic analysis of the orbits of this group action, which reveals a key geometric distribution. To describe this distribution globally, we introduce a canonical object, the **conformal co-symplectic structure** σ, defined by pushing forward the conformal class of the inverse of the original symplectic form. We prove that the image of this structure coincides with the geometric distribution identified previously. On the subspace of time-like and space-like geodesics, this structure is non-degenerate and defines a conformal class of symplectic forms. On the null subspace, its image is a codimension-1 distribution that we prove is the canonical contact structure on the space of light rays.

## INTRODUCTION

The space of geodesic curves on a Riemannian manifold possesses a canonical symplectic structure. On complete manifolds, this space of "parametrized geodesics" is isomorphic to the tangent bundle, and the physically meaningful space of "unparametrized geodesic trajectories" is then typically obtained via symplectic reduction. This elegant framework, however, confronts fundamental obstacles when one moves from the Riemannian to the pseudo-Riemannian setting.

*Date*: October 2, 2025.
2020 *Mathematics Subject Classification*. 53C50, 53D05, 53D10.
(P.I-Z thanks the Hebrew University of Jerusalem, Israel, for its continuous academic support. He is also grateful for the stimulating discussions and assistance provided by the AI assistant Gemini (Google).)

Let $(M, g)$ be a pseudo-Riemannian manifold. The space of parametrized geodesics, $\mathscr{G}_{\text{curv}}$, is still endowed with a canonical symplectic form $\omega$, and the dynamics are governed by the Hamiltonian $H(X, V) = \frac{1}{2} g(V, V)$. Unlike in the Riemannian case, its level sets partition $\mathscr{G}_{\text{curv}}$ into three dynamically distinct sectors: the space-like ($H > 0$), time-like ($H < 0$), and null ($H = 0$) geodesics. The standard method of symplectic reduction, which relies on quotienting the level sets of the Hamiltonian by its flow, consequently shatters the unified space of trajectories.

This fragmentation leads to a fatal obstruction. While the reduction process correctly produces symplectic structures on the spaces of space-like and time-like trajectories (of dimension $2n - 2$), it fails for the null case. The space of null trajectories has a natural dimension of $2n - 3$ but inherits by reduction a dimension of $2n - 2$, making it impossible to glue the resulting geometric spaces together in a coherent way. This critical failure leaves the geometry of null geodesics—essential for theories like General Relativity—disconnected from that of their massive counterparts, creating a puzzle that cannot be solved within the standard reduction framework.

To present our solution with maximal clarity, this paper will develop the argument in four-dimensional Minkowski spacetime. As a flat and complete manifold, its space of geodesics is analytically simple, allowing the core geometric construction to be demonstrated without distraction. This setting serves as the essential blueprint for the general theory, as the fundamental geometric principles we develop here are universal; the extension to curved and incomplete manifolds introduces significant analytical complexities but does not alter the core algebraic structure of the solution. The failure of the standard reduction method thus forces a change in strategy.

In the pseudo-Riemannian context, the physically meaningful space is not that of **geodesic curves** (or parametrized geodesics), but of **geodesic trajectories** (or unparametrized geodesic), as each trajectory represents a worldline with its own intrinsic proper time. Two parametrized geodesics, $\gamma(t) = X + tV$ and $\gamma'(t) = X' + tV'$, describe the same trajectory if and only if one is an affine reparametrization of the other. This means there must exist constants $a > 0$ and $b \in \mathbf{R}$ such that $\gamma'(t) = \gamma(at + b)$, which implies $X' = X + bV$ and $V' = aV$. The space of geodesic trajectories, $\mathscr{G}_{\text{traj}}$, is therefore constructed directly as the quotient of the space of parametrized geodesics by the action of the positive affine group $\text{Aff}^+(\mathbf{R})$:

$$\mathscr{G}_{\text{traj}} \cong \mathscr{G}_{\text{curv}} / \text{Aff}^+(\mathbf{R}).$$

This approach unifies all geodesic types but requires us to abandon the hope of finding a simple symplectic structure on the quotient, as the symplectic form $\omega$ is not invariant under affine reparametrizations; it is merely scaled.

This conformal scaling, while a failure from a purely symplectic viewpoint, can be overcome and lead to a unified solution. While $\omega$ and its inverse $\omega^{-1}$ are not invariant, the **conformal class** of the inverse, $[\omega^{-1}]$, *is* a perfect invariant of the action. This is the geometric object that survives the quotient. It can be pushed forward to $\mathscr{G}_{\text{traj}}$, defining the canonical **conformal co-symplectic structure**, $\sigma$.[1]

The central result of this paper is that this single structure $\sigma$ provides a unified geometric description of all geodesics. We will prove that on the open submanifolds of space-like and time-like trajectories, $\sigma$ is non-degenerate and defines a canonical conformal class of symplectic forms, while on the submanifold of null trajectories, its **image** is a codimension-1 distribution that is precisely the canonical **contact structure**.

**Historical Note.** The construction presented in this paper originates from the author's unpublished work, first formulated in a letter to Sergei Tabachnikov in May 2006. This private communication was subsequently cited in the literature, notably by B. Khesin and S. Tabachnikov in their 2009 paper on pseudo-Riemannian geodesics, and has since been incorporated into the author's forthcoming book, *The Geometry of Motion* (Appendix G). The present article is intended to provide the first formal, self-contained publication of this method, establishing the rigorous mathematical foundation for the conceptual argument presented in the book's appendices. A scanned copy of the original 2006 letter is available for reference [PIZ06].

## I. THE CANONICAL SYMPLECTIC STRUCTURE AND THE INADEQUACY OF REDUCTION

Before specializing to Minkowski spacetime, we briefly recall the foundational definitions of parametrized geodesics on a general Riemannian (or pseudo-Riemannian) manifold $(M, g)$. From a purely differential geometric perspective, a geodesic is a curve whose tangent vector is parallel-transported along itself, satisfying the condition of zero covariant acceleration.

A second, equivalent approach, which we adopt for its structural power, derives the geometry from a presymplectic framework. One considers the space of initial conditions $Y = TM \times \mathbf{R}$ and defines on it the Cartan 1-form:[2]

$$\varpi_y(\delta y) = g(v, \delta x) - \frac{1}{2} g(v, v) \delta t.$$

---

[1] The conformal class of the covariant tensor $\omega$ is also invariant under pullback. The difference, however, is that a structure on the quotient space must be defined via pushforward, an operation natural only to contravariant tensors. This makes the class $[\omega^{-1}]$ the only candidate.

[2] We use the notation $\delta y, \delta' y$, etc., to denote tangent vectors at the point $y$.

Its exterior derivative, $d\varpi$, is a presymplectic 2-form on Y. The parametrized geodesics are precisely the integral curves of the 1-dimensional characteristic distribution of $\ker d\varpi$. The space of parametrized geodesics, $\mathscr{G}_{\mathrm{curv}}$, is then defined as the quotient of Y by this characteristic foliation. By construction, this quotient space $\mathscr{G}_{\mathrm{curv}}$ is endowed with a unique symplectic form $\omega$ satisfying $\pi^*(\omega) = d\varpi$, where $\pi : Y \to \mathscr{G}_{\mathrm{curv}}$ is the quotient map.

This presymplectic approach has the advantage of immediately and canonically endowing the space of geodesic curves with its symplectic structure. However, for a general manifold, the geodesic flow may not be complete, in which case $\mathscr{G}_{\mathrm{curv}}$ may fail to be a Hausdorff manifold. This is a primary motivation for focusing our analysis on Minkowski spacetime, where the metric is complete, the geodesic flow is global, and $\mathscr{G}_{\mathrm{curv}}$ is a well-behaved manifold that can be identified with TM−null section, equipped with its standard symplectic form. This allows us to isolate the fundamental geometric problem of the quotient structure from these secondary analytical complexities.

**1. The Space of Parametrized Geodesics.** In Minkowski spacetime $(\mathrm{M}, g) \cong (\mathbf{R}^4, g)$, where $g$ has signature $(+, +, +, -)$, a geodesic is a parametrized affine straight line, $\gamma(t) = \mathrm{X} + t\mathrm{V}$. The space of all such parametrized geodesics is the smooth, 8-dimensional manifold:

$$\mathscr{G}_{\mathrm{curv}} \cong \mathrm{TM} - \text{null section} = \mathbf{R}^4 \times (\mathbf{R}^4 - \{0\}).$$

**2. The Canonical Symplectic Form.** The space $\mathscr{G}_{\mathrm{curv}}$ is endowed with a canonical symplectic structure. Let $\gamma = (\mathrm{X}, \mathrm{V}) \in \mathscr{G}_{\mathrm{curv}}$. The canonical symplectic form $\omega$ on $\mathscr{G}_{\mathrm{curv}}$ is given by:

$$\omega_\gamma(\delta\gamma, \delta'\gamma) = g(\delta\mathrm{V}, \delta'\mathrm{X}) - g(\delta'\mathrm{V}, \delta\mathrm{X}).$$

This 2-form is non-degenerate and exact, with $\omega = d\varpi$, where $\varpi$ is the Liouville 1-form defined by $\varpi_\gamma(\delta\gamma) = g(\mathrm{V}, \delta\mathrm{X})$. The pair $(\mathscr{G}_{\mathrm{curv}}, \omega)$ is therefore a symplectic manifold.

**3. The Obstruction from Hamiltonian Reduction.** The standard method for passing from parametrized curves to unparametrized trajectories is Hamiltonian reduction. This procedure breaks down in the pseudo-Riemannian context. The Hamiltonian $\mathrm{H} = \frac{1}{2}g(\mathrm{V}, \mathrm{V})$ is indefinite, partitioning $\mathscr{G}_{\mathrm{curv}}$ into three disjoint regions. Applying reduction level-set by level-set produces three disconnected symplectic manifolds of dimension 6.

A fundamental inconsistency arises for the null geodesics. The space of all affine lines in $\mathbf{R}^4$ is a 6-dimensional manifold. Within this space, the space of *null* lines is a 5-dimensional submanifold. The reduction procedure, however, yields a 6-dimensional symplectic manifold for the null case.[3] This creates an insurmountable topological obstruction, as a 6-dimensional space cannot serve as the 5-dimensional boundary of the other two spaces. This dimensional conflict forces a choice: either accept the fragmentation, or seek a new method.

## II. Symplectic Analysis of the Reparametrization Orbits

To understand the geometry on the quotient $\mathscr{G}_{\text{traj}} = \mathscr{G}_{\text{curv}}/\text{Aff}^+(\mathbf{R})$, we must first analyze the symplectic nature of the partition of $\mathscr{G}_{\text{curv}}$ into the orbits of the affine group. This analysis is the key to the solution.

Note first that the action of $\text{Aff}^+(\mathbf{R})$ on $\mathscr{G}_{\text{curv}}$ is free and all the orbits are 2-dimensional. Remark also that quotient map

$$\pi : \mathscr{G}_{\text{curv}} \to \mathscr{G}_{\text{traj}}$$

is then a principal fiber bundle, and since $\text{Aff}^+(\mathbf{R})$ is contractible, a trivial $\text{Aff}^+(\mathbf{R})$-principal bundle.

Let $\gamma \in \mathscr{G}_{\text{curv}}$. The orbit map is $\hat{\gamma} : (a,b) \mapsto (a,b) \cdot \gamma$.

**Proposition 1.** *The pullback of the symplectic form $\omega$ by the orbit map is proportional to the standard symplectic form on $\text{Aff}^+(\mathbf{R})$:*

$$\hat{\gamma}^*(\omega) = g(V,V)\, da \wedge db.$$

*Proof.* The orbit map is $\hat{\gamma}(a,b) = (X + bV, aV)$. The tangent vectors to the orbit corresponding to the coordinates $(a,b)$ are generated by the partial derivatives:

$$\xi_a = \frac{\partial \hat{\gamma}}{\partial a} = (0, V) \quad \text{and} \quad \xi_b = \frac{\partial \hat{\gamma}}{\partial b} = (V, 0).$$

The pullback form acts on these tangent vectors as $\hat{\gamma}^*(\omega)(\xi_a, \xi_b) = \omega_{\hat{\gamma}(a,b)}(\xi_a, \xi_b)$. Applying the definition of $\omega$:

$$\omega((0,V),(V,0)) = g(V,V) - g(0,0) = g(V,V).$$

Since $(da \wedge db)(\xi_a, \xi_b) = 1$, we obtain the identity. $\square$

This has immediate consequences for the geometry of the orbits:

- For **non-null** geodesics, $g(V,V) \neq 0$. The orbits $\mathscr{O}_\gamma$ are 2-dimensional symplectic submanifolds of $\mathscr{G}_{\text{curv}}$.
- For **null** geodesics, $g(V,V) = 0$. The orbits are **isotropic**, meaning the symplectic form vanishes identically when restricted to them.

---

[3]The symplectic reduction adds a dimension, sometimes called "color," to the null trajectories. While this parameter may have a role in quantization, it is superfluous in this purely geometric setting.

This distinction motivates the study of the distribution of symplectic orthogonal spaces to the orbits. Let

$$F_\gamma = \text{Orth}_\omega(T_\gamma \mathcal{O}_\gamma)$$

be the symplectic orthogonal to the tangent space of the orbit at $\gamma$. We define the distribution $\mathscr{F}$ on $\mathscr{G}_{\text{traj}}$ as the pushforward of F:

$$\mathscr{F}_\tau = D\pi_\gamma(F_\gamma), \quad \text{where } \tau = \pi(\gamma).$$

The nature of this distribution depends directly on the orbit type.

- In the non-null case, the orbit is symplectic, which implies the direct sum decomposition $T_\gamma \mathscr{G}_{\text{curv}} = T_\gamma \mathcal{O}_\gamma \oplus F_\gamma$. The pushforward $D\pi_\gamma$ maps $F_\gamma$ isomorphically onto the tangent space of the quotient, so $\mathscr{F}_\tau = T_\tau \mathscr{G}_{\text{traj}}$.
- In the null case, the orbit is isotropic, so $T_\gamma \mathcal{O}_\gamma \subset F_\gamma$. The distribution $F_\gamma$ is co-isotropic. The quotient of a co-isotropic space by its isotropic kernel is a symplectic space, but here we quotient by the full orbit. The resulting distribution $\mathscr{F}_\tau$ is a proper subspace of codimension 1.

This analysis thus reveals the fundamental geometric structures inherited by the quotient: a full-rank symplectic distribution on the space of massive trajectories $\mathscr{G}_{\text{traj}}^\pm$, and a codimension-1 distribution on the space of null trajectories $\mathscr{G}_{\text{traj}}^0$. The task is now to find a single, canonical object on $\mathscr{G}_{\text{traj}}$ that unifies these two distinct geometric behaviors.

### III. UNIFICATION VIA THE CONFORMAL CO-SYMPLECTIC STRUCTURE

The fragmented picture provided by the symplectic analysis can be unified by considering the inverse of the symplectic form, the contravariant tensor $\omega^{-1}$.[4]

**4. The Canonical Conformal Co-symplectic Structure.** A direct calculation shows that $\omega^{-1}$ is not invariant under the action of $\text{Aff}^+(\mathbf{R})$, but transforms conformally.

**Proposition 2.** *For any* $(a, b) \in \text{Aff}^+(\mathbf{R})$*:*

$$(a, b)_*(\omega^{-1}) = a \cdot \omega^{-1}.$$

---

[4]This tensor has deep historical roots. Its components are precisely the **Lagrange parentheses**, introduced by Lagrange in his seminal 1808 work [Lag08] on the variation of constants —a paper that marks the genesis of symplectic geometry. Our method is thus, in a deep sense, a return to the very first geometric object uncovered by analytical mechanics, an object whose contravariant nature proves essential for the pushforward construction required by the quotient. For a detailed discussion of the historical priority of Lagrange in this context, see *Lagrange et Poisson, sur la Variation des Constantes* by the author in *Siméon-Denis Poisson, Les Mathématiques au service de la science*, Yvette Kosmann-Schwarzbach (éd.), Palaiseau, Éditions de l'École polytechnique, 2013.

*Proof.* Notice first that, for all $(a,b) \in \text{Aff}^+(\mathbf{R})$:

$$(a,b)^*(\omega) = a\omega.$$

Indeed, for two tangent vectors $\delta\gamma$ and $\delta'\gamma$, we have:

$$(a,b)^*(\omega)(\delta\gamma, \delta'\gamma) = g(\delta[aV], \delta'[X+bV]) - g(\delta'[aV], \delta[X+bV])$$
$$= g(a\delta V, \delta' X) + abg(\delta V, \delta' V) - g(a\delta' V, \delta X) - abg(\delta' V, \delta V)$$
$$= ag(\delta V, \delta' X) - ag(\delta' V, \delta X) = a\omega(\delta\gamma, \delta'\gamma)$$

Now, let's use matrix notation to simplify the computation. Let $\omega(\xi, \xi') = \xi^t \Omega \xi'$, where tangent vectors are column matrices and the superscript $t$ denotes transposition. The inverse tensor $\omega^{-1}$ is represented by $\Omega^{-1}$. If 1-forms are represented by row matrices, the evaluation $\omega^{-1}(\alpha, \beta) = \beta(\omega^{-1}(\alpha))$ corresponds to the matrix product $\beta\Omega^{-1}\alpha^t$. Hence, for any diffeomorphism $\phi$, the pullback $\phi^*(\omega)$ is represented by the matrix $M^t \Omega M$, with $M = D\phi$ and $\phi_*(\omega^{-1})$ is represented by $M\Omega^{-1}M^t$. For $\phi = (a,b)$, $\phi^{-1} = (1/a, -b/a)$. Hence, $M^t\Omega M = a\Omega$ and $M\Omega^{-1}M^t = [(M^{-1})^t \Omega M^{-1}]^{-1} = [(1/a)\Omega]^{-1} = a\Omega^{-1}$. Therefore $(a,b)_*(\omega^{-1}) = a\omega^{-1}$. □

It follows from this proposition that the **conformal class**, $[\omega^{-1}]$, *is* invariant. Hence, this invariant class descends to the quotient and defines the fundamental geometric structure on the space of trajectories.

**Definition.** *The **canonical conformal co-symplectic structure**[5] on $\mathcal{G}_{traj}$ is the field of rays of contravariant, antisymmetric 2-tensors, denoted $\sigma$, obtained by the pushforward of the conformal class of $\omega^{-1}$ by the quotient map $\pi$:*

$$\sigma := \pi_*([\omega^{-1}]).$$

This structure provides the crucial link between the symplectic analysis and a unified object. As we requested, here is the proof of the identity between the distribution $\mathcal{F}$ and the image of $\sigma$.

**Proposition.** *The distribution $\mathcal{F}_\tau$ is the image of the structure $\sigma_\tau$. That is,*

$$\mathcal{F}_\tau = \text{Im}(\sigma_\tau).$$

*Proof.* By definition, $\mathcal{F}_\tau = D\pi_\gamma(\text{Orth}_\omega(T_\gamma \mathcal{O}_\gamma))$. The image of $\sigma_\tau$ is the pushforward of the image of any representative of $[\omega^{-1}]$ acting on covectors that are pullbacks from the base: $\text{Im}(\sigma_\tau) = D\pi_\gamma(\text{Im}(\omega^{-1} \circ \pi^*))$. The space of pulled-back covectors, $\pi^*(T^*_\tau \mathcal{G}_{traj})$, is precisely the annihilator of the kernel of the differential, $(T_\gamma \mathcal{O}_\gamma)^\circ$. In symplectic linear algebra, for any subspace W, the image of its annihilator under the co-symplectic map is its symplectic orthogonal: $\omega^{-1}(W^\circ) = \mathrm{Orth}_\omega(W)$. Applying this, we have:

$$\mathrm{Im}(\omega^{-1} \circ \pi^*) = \omega^{-1}((T_\gamma \mathcal{O}_\gamma)^\circ) = \mathrm{Orth}_\omega(T_\gamma \mathcal{O}_\gamma) = F_\gamma.$$

---

[5]We have chosen this terminology to be precise. It is "co-symplectic" because it is a contravariant tensor, the inverse of a symplectic form. It is "conformal" because the tensor itself transforms by a scaling factor, meaning its conformal class is the true invariant. This should not be confused with the more common notion of a locally conformal symplectic (LCS) structure, which is a covariant property of a 2-form $\omega$ satisfying $d\omega = \theta \wedge \omega$ for some 1-form $\theta$.

Pushing this identity forward via $D\pi_\gamma$ yields: $\mathrm{Im}(\sigma_\tau) = D\pi_\gamma(F_\gamma) = \mathscr{F}_\tau$. □

**5. The Main Theorem.** With the identity $\mathscr{F} = \mathrm{Im}(\sigma)$ established, we can now state our main theorem, which unifies the geometry of all geodesic trajectories.

**Theorem.** *The canonical conformal co-symplectic structure $\sigma$ on $\mathscr{G}_{traj}$ has the following properties:*

(1) *On the open submanifolds of space-like and time-like trajectories $\mathscr{G}^\pm_{traj}$, its image is the full tangent space, $\mathrm{Im}(\sigma) = T\mathscr{G}^\pm_{traj}$. The structure $\sigma$ is non-degenerate and defines a canonical conformal class of symplectic forms.*
(2) *On the submanifold of null trajectories $\mathscr{G}^0_{traj}$, its image $\mathrm{Im}(\sigma)$ is a co-dimension-1 distribution that is the canonical contact structure.*

*Proof.* Part (1) follows directly from the symplectic analysis in Section II. For any $\tau \in \mathscr{G}^\pm_{traj}$, the distribution $\mathscr{F}_\tau$ is the full tangent space, $T_\tau \mathscr{G}_{traj}$. Since $\mathscr{F}_\tau = \mathrm{Im}(\sigma_\tau)$, the structure $\sigma$ is non-degenerate on these open submanifolds. Any representative tensor from the class $\sigma_\tau$ is therefore a co-symplectic tensor, and its inverse is a symplectic form, endowing $\mathscr{G}^\pm_{traj}$ with a canonical conformal class of symplectic structures.

For Part (2), we must prove that the distribution $\mathscr{F}_\tau = \mathrm{Im}(\sigma_\tau)$ on $\mathscr{G}^0_{traj}$ is contact. To do so, we use a global cross-section $\Sigma$ of the bundle $\pi : \mathscr{G}^0_{curv} \to \mathscr{G}^0_{traj}$, the existence of which is guaranteed by the simple topology of Minkowski spacetime, though its construction presents a significant analytical challenge in a general curved spacetime. The distribution $\mathscr{F}_\tau$ on $\mathscr{G}^0_{traj}$ corresponds to the distribution $D_F = F_\gamma \cap T_\gamma \Sigma$ on the section. We will prove that $D_F$ is the kernel of a contact 1-form.

Let $\alpha = \varpi|_\Sigma$ be the restriction of the Liouville 1-form to the section. A tangent vector $\eta = (\delta X, \delta V)$ belongs to $D_F$ if it is tangent to $\Sigma$ and symplectically orthogonal to the orbit generators $\xi_T = (V, 0)$ and $\xi_S = (0, V)$. Orthogonality to $\xi_T$ is automatically satisfied for any vector tangent to $\Sigma$. Orthogonality to $\xi_S$ requires $\omega(\eta, \xi_S) = -g(V, \delta X) = 0$. This condition is precisely the definition of the kernel of the 1-form $\alpha(\eta) = g(V, \delta X)$. Therefore, the two distributions coincide: $D_F = \ker(\alpha)$.

The final step is to show that $\alpha$ is a contact form. The characteristic distribution of $\omega = d\varpi$ on $\mathscr{G}^0_{curv}$ is spanned by $\xi_T = (V, 0)$. Since the section $\Sigma$ is everywhere transverse to the orbits, $\xi_T$ is nowhere tangent to $\Sigma$. This implies that the restriction of the characteristic distribution to the section is trivial:

$\ker(\varpi|_\Sigma) \cap \ker(d\varpi|_\Sigma) = \{0\}$. This is the defining condition for $\alpha = \varpi|_\Sigma$ to be a contact form.[6] □

## CONCLUSION

The inadequacy of the standard reduction method was not a flaw, but a signal of a richer underlying geometry. By first analyzing the symplectic geometry of the reparametrization orbits, we identified a crucial distribution, $\mathscr{F}$. We then constructed the canonical conformal co-symplectic structure $\sigma$ and proved that $\mathscr{F}$ is its image. This single, unified object $\sigma$ behaves exactly as required across the entire space of trajectories: on the open domains of space-like and time-like geodesics, its image is the full tangent space, endowed with a canonical class of symplectic forms, while on the submanifold of null geodesics, its image is precisely the canonical contact distribution.

The payoff for this effort is a profound unification: the canonical conformal co-symplectic structure reveals the symplectic geometry of massive particles and the contact geometry of light as two faces of a single, underlying geometric truth. This unification is not merely a matter of geometric elegance; it offers a new perspective on the foundations of geometric quantization and representation theory. The submanifold of null trajectories, $\mathscr{G}^0_{\text{traj}}$, is precisely the contact manifold whose **symplectization** yields the 6-dimensional coadjoint orbit of the Poincaré group representing a massless particle. Viewed differently, $\mathscr{G}^0_{\text{traj}}$ serves as the natural **prequantum bundle** over the 4-dimensional symplectic manifold of unlocated light rays —the covariant "celestial sphere of spacetime."

The novelty of the present work lies in revealing that these structures are intrinsic consequences of pseudo-Riemannian geometry itself, not merely artifacts of Minkowski spacetime's symmetries. The algebraic framework of the conformal co-symplectic structure is general, and this generality is crucial because it corrects a foundational modeling error in the standard approach. The standard method implicitly assumes that the tangent bundle TM provides a complete model for the space of all parametrized geodesics. This assumption, however, is flawed—an issue not unique to the pseudo-Riemannian setting, but one that arises even in simple Riemannian manifolds such as the plane with a point removed. The space of geodesics in such cases is not Hausdorff, and the tangent bundle, as a single slice of initial data, fails to capture the totality of geodesic trajectories.

This fundamental inadequacy is severely amplified in the context of incomplete pseudo-Riemannian manifolds, where the problem of null geodesics adds to the difficulty. In this situation, the standard theorems of global symplectic geometry become inapplicable, and the reduction of the tangent bundle loses its physical meaning, ignoring rather than solving the core topological pathologies. Our approach, based on the algebraic action of affine reparametrizations on the space of *all* geodesics, is intrinsically suited to this general setting. By operating at the level of germs of geodesics, it bypasses these global and modeling obstructions from the outset. The Minkowski spacetime model thus serves as the essential, non-trivial blueprint for this robust and general theory. This perspective opens a promising avenue for understanding how these structures might generalize to curved spacetimes, suggesting that quantization is not an additional layer imposed upon mechanics, but a structure already woven into the geometric fabric of spacetime trajectories. These implications will be explored in more detail in the author's forthcoming book, *The Geometry of Motion*.

---

[6]Which is equivalent to the non-degeneracy requirement that the top-degree form $\alpha \wedge (d\alpha)^k$ defines a volume form on the $(2k+1)$-dimensional manifold.

P. IGLESIAS-ZEMMOUR, EINSTEIN INSTITUTE OF MATHEMATICS, THE HEBREW UNIVERSITY OF JERUSALEM, EDMOND J. SAFRA CAMPUS, GIVAT RAM, 9190401 JERUSALEM, ISRAEL.

*Email address*: piz@math.huji.ac.il

*URL*: http://math.huji.ac.il/~piz